\documentclass[12pt,a4paper]{article}
\usepackage{amsfonts}
\usepackage{amsmath}

\newtheorem{theorem}{Theorem}[section]

\newtheorem{prop}{Proposition}[section]
\newtheorem{cor}{Corollary}[section]

\newcommand{\R}{\mathbb{R}}

\newcommand{\nN}{n \in \mathbb{N}}

\newcommand{\C}{\mathbb{C}}

\newcommand{\g}{\textbf{g}}

\newcommand{\tr}{\mbox{tr}}

\newcommand{\bean}{\begin{eqnarray*}}
\newcommand{\eean}{\end{eqnarray*}}

\newcommand{\la}{\langle}
\newcommand{\ra}{\rangle}

\newcommand{\G}{\widehat{G}}

\newcommand{\di}{\mbox{\mbox{dim}}}

\date{}

\begin{document}

\title{Smoothness of Densities on Compact Lie Groups \footnote{For the proceedings of the 8th Isaac Congress, Moscow 2011.}}

\author{ David Applebaum, \\ School of Mathematics and Statistics,\\ University of
Sheffield,\\ Hicks Building, Hounsfield Road,\\ Sheffield,
England, S3 7RH\\ ~~~~~~~\\e-mail: D.Applebaum@sheffield.ac.uk }

\maketitle
\begin{abstract}
We give necessary and sufficient conditions for both square integrability and smoothness for densities of a probability measure on a compact connected Lie group.
\end{abstract}

\section{Introduction}

The study of probability measures on groups provides a mathematical framework for describing the interaction of chance with symmetry. This subject is broad and interacts with many other areas of mathematics and its applications such as analysis on groups \cite{Var}, stochastic differential geometry \cite{El}, statistics \cite{Dia} and engineering \cite{Chir}.

In this paper we focus on the important question concerning when a probability measure on a compact group has a regular density with respect to Haar measure. We begin by reviewing work from \cite{App1} where Peter-Weyl theory is used to find a necessary and sufficient condition for such a measure to have a square-integrable density. This condition requires the convergence of an infinite series of terms that are formed from the (non-commutative) Fourier transform of the measure in question. We also describe a related result from \cite{App2} where it is shown that square-integrability of the measure is a necessary and sufficient condition for the associated convolution operator to be Hilbert-Schmidt (and hence compact) on the $L^{2}$-space of Haar measure.

In the second part of our paper we turn our attention to measures with smooth densities. A key element of our approach is the important insight of Hermann Weyl that the unitary dual $\G$ of the group $G$ can be parameterised by the space of highest weights. This effectively opens up $\G$ to investigation by standard analytical methods. We introduce Suguira's space of rapidly decreasing functions of weights which was shown in \cite{Sug} to be topologically isomorphic to $C^{\infty}(G)$. We are then able to prove that a probability measure has a smooth density if and only if its Fourier transform lives in Suguira's space. This improves on results of \cite{App3} where the Sobolev embedding theorem was used to find sufficient conditions for such a density to exist.

In the last part of the paper we give a brief application to statistical inference. In \cite{KR}, Kim and Richards have introduced an estimator for the density of a signal on the group based on i.i.d. (i.e. independent and identically distributed) observations of the signal after it has interacted with an independent noise. To obtain fast rates of convergence to the true density, the noise should be in a suitable ``smoothness class'' where smoothness is here measured in terms of the decay of the Fourier transform of the measure. We show that the ``super-smooth'' class is smooth in the usual mathematical sense.

\section{Fourier Transforms of Measures on Groups }

Throughout this paper $G$ is a compact connected Lie group with neutral element $e$ and dimension $d$, ${\cal B}(G)$ is the Borel $\sigma$-algebra of $G$ and ${\cal P}(G)$ is the space of probability measures on $(G, {\cal B}(G))$, equipped with the topology of weak convergence. The role of the uniform distribution on $G$ is played by {\it normalised Haar measure} $m \in {\cal
P}(G)$ and we recall that this is a bi--invariant measure in that
$$ m(A\sigma) = m(\sigma A) = m(A),$$
for all $A \in {\cal B}(G), \sigma \in G$. We will generally write $m(d\sigma) = d\sigma$ within integrals.

Our main focus in this paper is those $\rho \in {\cal P}(G)$ that are absolutely continuous with respect to $m$ and so they have densities $f \in L^{1}(G,m)$ satisfying
       $$ \rho(A) = \int_{A}f(\sigma)d\sigma, $$
       for all $A \in {\cal B}(G)$.









 A key tool which we will use to study these measures is the non-commutative Fourier transform which is defined using representation theory. We recall some key facts that we need. A good reference for the material below about group representations, the Peter-Weyl theorem and Fourier analysis of square-integrable functions is Faraut \cite{Far}.

If $H$ is a complex separable Hilbert space then ${\cal U}(H)$ is
the group of all unitary operators on $H$. A {\it unitary representation} of $G$ is a strongly continuous homomorphism $\pi$
from $G$ to ${\cal U}(V_{\pi})$ for some such Hilbert space $V_{\pi}$. So we have for all $g,h \in G,$:
\begin{itemize}
\item $ \pi(gh) = \pi(g)\pi(h),$
 \item $ \pi(e) = I_{\pi}$ (where $I_{\pi}$ is the identity operator on $V_{\pi}$,)
\item $\pi(g^{-1}) = \pi(g)^{-1} = \pi(g)^{*}.$
\end{itemize}

$\pi$ is {\it irreducible} if it has no non-trivial invariant closed subspace.
Every group has a trivial representation $\delta$ acting on $V_{\delta} = \C$ by $\delta(g) = 1$ for all
$g \in G$ and it is clearly irreducible. The {\it unitary dual} of $G, \widehat{G}$ is defined to be the set of equivalence classes of all irreducible
representations of $G$ with respect to unitary conjugation. We will as usual identify each equivalence class with a typical representative element. As $G$ is compact, for all $\pi \in \G, d_{\pi}:=\di(V_{\pi}) < \infty$ so that each $\pi(g)$ is a unitary matrix. Furthermore in this case $\widehat{G}$ is countable.

For each $\pi \in \G$, we define co-ordinate functions $\pi_{ij}(\sigma) =
\pi(\sigma)_{ij}$ with respect to a some orthonormal basis in $V_{\pi}$.

\begin{theorem} [Peter-Weyl] \label{PW} The set $\{\sqrt{d_{\pi}}\pi_{ij}, 1 \leq i,j \leq
d_{\pi}, \pi \in \G\}$ is a complete orthonormal basis for
$L^{2}(G,\C)$.
\end{theorem}

The following consequences of Theorem \ref{PW} are straightforward to derive using Hilbert space arguments.

\begin{cor} \label{Pyth} For $f,g \in L^{2}(G,\C)$

\begin{itemize} \item Fourier expansion.

$$ f = \sum_{\pi \in \G}d_{\pi}\tr(\widehat{f}(\pi)\pi),$$

where $\widehat{f}(\pi): =
\int_{G}f(\sigma^{-1})\pi(\sigma)d\sigma$ is the Fourier transform of $f$.

\item The Plancherel theorem.
$$ ||f||^{2} = \sum_{\pi \in
\G}d_{\pi}|||\widehat{f}(\pi)|||^{2}$$

where $|||\cdot|||$ is the Hilbert-Schmidt norm $|||T|||: =
\tr(TT^{*})^{\frac{1}{2}}$.

\item The Parseval identity.
$$ \la f, g \ra = \sum_{\pi \in
\G}d_{\pi}\tr(\widehat{f}(\pi)\widehat{g}(\pi)^{*}).$$

\end{itemize}

\end{cor}

If $\mu \in {\cal P}(G)$ we define its {\it Fourier transform} $\widehat{\mu}$ to be
$$ \widehat{\mu}(\pi) = \int_{G}\pi(\sigma^{-1})\mu(d\sigma),$$
for each $\pi \in \G$. For example if $\epsilon_{e}$ is a Dirac mass at $e$ then $\widehat{\epsilon_{e}}(\pi) = I_{\pi}$ and $\widehat{m}(\pi) = \left\{\begin{array}{c c}
                          0 & \mbox{if}~\pi \neq \delta\\
                          1 & \mbox{if}~\pi = \delta
                          \end{array}\right.$.
If $\mu$ has a density $f$ then $\widehat{\mu} = \widehat{f}$ as defined in Corollary \ref{Pyth}. If we take $G$ to be the $d$-torus $\mathbb{T}^{d}$ then $\widehat{G}$ is the dual group $\mathbb{Z}^{d}$ and the Fourier transform is precisely the usual characteristic function of the measure $\mu$ defined by $\widehat{\mu}(n) = \int_{\mathbb{T}^{d}}e^{-i n \cdot x}\mu(dx)$ for $n \in \mathbb{Z}^{d}$, where $\cdot$ is the scalar product. Note that any compact connected abelian Lie group is isomorphic to $\mathbb{T}^{d}$.

Fourier transforms of measures on groups have been studied by many authors, see e.g. \cite{KI, He1, He2, Sieb} where proofs of the following basic properties can be found.
\vspace{5pt}

For all $\mu, \mu_{1}, \mu_{2} \in {\cal P}(G), \pi \in \G$,
 \begin{enumerate} \item  $\widehat{\mu_{1} *
\mu_{2}}(\pi) = \widehat{\mu_{2}}(\pi)\widehat{\mu_{1}}(\pi),$
\item $\widehat{\mu}$ determines $\mu$ uniquely, \item
$||\widehat{\mu}(\pi)||_{\infty} \leq 1$, where $||\cdot||_{\infty}$ denotes the operator norm in $V_{\pi}$. \item Let $(\mu_{n}, \nN)$ be a sequence in ${\cal P}(G)$. $\mu_{n}
\rightarrow \mu$ (weakly) if and only if $\widehat{\mu_{n}}(\pi)
\rightarrow \widehat{\mu}(\pi)$. \end{enumerate}

\vspace{5pt}

{\bf Remark.} Most authors define $\widehat{\mu}(\pi) = \int_{G}\pi(\sigma)\mu(d\sigma)$. This has the advantage that Property 1 above will then read  $\widehat{\mu_{1} * \mu_{2}}(\pi) = \widehat{\mu_{1}}(\pi)\widehat{\mu_{2}}(\pi)$ but the disadvantage that if $\mu$ has density $f$ then $\widehat{\mu}(\pi) = \widehat{f}(\pi)^{*}$. It is also worth pointing out that the Fourier transform continues to make sense and is a valuable probabilistic tool in the case where $G$ is a general locally compact group (see e.g. \cite{He1, He2, Sieb}.)

\section{Measures With Square-Integrable Densities}

In this section we examine the case where $\mu$ has a square-integrable density. The following result can be found in \cite{App1} and so we only sketch the proof here.

\begin{theorem} \label{sqrint} $\mu$ has an $L^{2}$-density $f$ if and only if
  $$ \sum_{\pi \in \G}d_{\pi}|||\widehat{\mu}(\pi)|||^{2} < \infty.$$ In this case
  $$ f(\sigma) = \sum_{\pi \in \G}d_{\pi}\tr(\widehat{\mu}(\pi)\pi(\sigma))~\mbox{a.e.}.$$
\end{theorem}

{\it Proof.} Necessity is straightforward. For sufficiency define $g: = \sum_{\pi \in
\G}d_{\pi}\tr(\widehat{\mu}(\pi)\pi).$ Then $g \in L^{2}(G, \C)$ and by uniqueness of Fourier
coefficients $\widehat{g}(\pi) = \widehat{\mu}(\pi)$. Using
Parseval's identity, Fubini's theorem and Fourier expansion, we find that for each $h \in C(G, \C)$:
$$ \int_{G}h(\sigma)\overline{g(\sigma)}d\sigma = \sum_{\pi \in
\G}d_{\pi}\tr(\widehat{h}(\pi)\widehat{\mu}(\pi)^{*}) =
\int_{G}h(\sigma)\mu(d\sigma).$$ This together with the Riesz
representation theorem implies that $g$ is real valued and
$g(\sigma)d\sigma = \mu(d\sigma)$. The fact that $g$ is non-negative then follows from the Jordan decomposition for signed measures.
$\hfill \Box$

\vspace{5pt}

See \cite{App1} for specific examples. We will examine some of these in the next section from the finer point of view of smoothness.

To study random walks and L\'{e}vy processes in $G$ we need the convolution
operator $T_{\mu}$ in $L^{2}(G,\C)$ associated to $\mu \in {\cal P}(G)$ by
$$ (T_{\mu}f)(\sigma): = \int_{G}f(\sigma \tau) \mu(d\tau),$$
for $f \in L^{2}(G, \C), \sigma \in G$. For example $T_{\mu}$ is the transition operator corresponding to the random walk $(\mu^{*n}, \nN)$. The following properties are fairly easy to establish.

\begin{itemize}
\item $T_{\mu}$ is a contraction. \item $T_{\mu}$ is self-adjoint
if and only if $\mu$ is symmetric, i.e. $\mu(A) = \mu(A^{-1})$ for all $A \in {\cal B}(G)$.
\end{itemize}

The next result is established in \cite{App2}.

\begin{theorem} \label {HS} The operator $T_{\mu}$ is Hilbert-Schmidt if and only if $\mu$ has a
square-integrable density.
\end{theorem}

{\it Proof.} Sufficiency is obvious by the Hilbert-Schmidt theorem.
For necessity, suppose that $T_{\mu}$ Hilbert-Schmidt. Then it has a kernel $k \in L^{2}(G \times G)$ and
   $$ (T_{\mu}f)(\sigma) = \int_{G}f(\sigma)k_{\mu}(\sigma,
   \tau)d\tau.$$
In particular for each $A \in {\cal B}(G)$,
$$ \mu(A) = T_{\mu}1_{A}(e) = \int_{A}k_{\mu}(e, \tau)d\tau.$$
It follows that $\mu$ is absolutely continuous with respect to $m$ with density $f= k_{\mu}(e, \cdot). \hfill
\Box$

\vspace{5pt}

Let $(\mu_{t}, t \geq 0)$ be a weakly continuous convolution
semigroup in ${\cal P}(G)$ and write $T_{t}: = T_{\mu_{t}}$. Then
$(T_{t}, t \geq 0)$ is a strongly continuous contraction semigroup
on $L^{2}(G, \C)$ (see e.g. \cite{Hu, He1, Liao, App2}.)

\begin{cor} The linear operator \label{TC} $T_{t}$ is trace-class for all $t > 0$ if and only if
$\mu_{t}$ has a square-integrable density for all $t > 0$.
\end{cor}

{\it Proof}. For each $t > 0$, if $\mu_{t}$ has a square-integrable density then $T_{t} = T_{\frac{t}{2}}T_{\frac{t}{2}}$ is the product of
two Hilbert-Schmidt operators and hence is trace class. The converse follows from the fact that every trace-class operator is Hilbert-Schmidt. $\hfill \Box$

\vspace{5pt}

If for $t > 0, \mu_{t}$ has a square-integrable density and is symmetric, then by Theorem \ref{HS}, $T_{t}$  is a compact self-adjoint operator and so has a discrete spectrum of positive eigenvalues
$1 = e^{-t\beta_{1}} > e^{-t\beta_{2}} > \cdots > e^{-t\beta_{n}} \rightarrow 0$ as $n \rightarrow \infty$. Furthermore by Corollary \ref{TC}, $T_{t}$ is trace class and $$\mbox{Tr}(T_{t}) = \sum_{n=1}^{\infty}e^{-t\beta_{n}} < \infty.$$ Further consequences of these facts including the application to small time asymptotics of densities can be found in \cite{App2, App3}.

\section{Sugiura Space and Smoothness}





In this section we will review key results due to Sugiura \cite{Sug} which we will apply to densities in the next section. In order to do this we need to know about weights on Lie algebras and we will briefly review the necessary theory.

\subsection{Weights}

Let $\g$ be the Lie algebra of $G$ and $\exp: \g \rightarrow G$ be
the exponential map. For each unitary representation $\pi$ of $G$ we obtain a Lie algebra
representation $d\pi$ by
$$ \pi(\exp(tX)) = e^{td\pi(X)}~\mbox{for all}~t \in \R.$$
Each $d\pi(X)$ is a skew-adjoint matrix on $V_{\pi}$ and
$$ d\pi([X,Y]) = [d\pi(X), d\pi(Y)],$$ for all $X,Y \in \g$.
A maximal torus $\mathbb{T}$ in $G$ is a maximal commutative
subgroup of $G$. Its dimension $r$ is called the rank of $G$. Here are some key facts about maximal tori.
\begin{itemize} \item Any $\sigma \in G$ lies on some maximal torus.
\item Any two maximal tori are conjugate.
\end{itemize}
Let $\mathrm{t} $ be the Lie algebra of $\mathbb{T}$. Then it is a maximal abelian subalgebra of $\g$. The matrices $\{d\pi(X), X \in \mathrm{t}\}$ are mutually commuting and so simultaneously diagonalisable, i.e. there exists a non-singular matrix $Q$ such that
$$ Q d \pi(X)Q^{-1} = \mbox{diag}(i\lambda_{1}(X), \ldots, i\lambda_{d_{\pi}}(X)).$$ The distinct linear functionals $\lambda_{j}$ are called the {\it weights} of $\pi$.

Let $Ad$ be the adjoint representation of $G$ on $\g$.  We can and will choose an $Ad$-invariant inner product $(\cdot, \cdot)$ on $\g$. This induces an inner
product on $\mathrm{t}^{*}$ the algebraic dual of $\mathrm{t}$ which we also write as $(\cdot, \cdot)$. We denote the corresponding norm by $|\cdot|$.
The weights of the adjoint representation acting on $\g$ equipped with $(\cdot, \cdot)$ are called the {\it roots} of $G$. Let ${\mathcal P}$ be the set of all roots of $G$. We choose a convention for positivity of roots as follows. Pick $v \in \mathrm{t}$ such that ${\mathcal P} \cap \{\eta \in \mathrm{t}^{*}; \eta(v)\} = \emptyset$. Now define ${\mathcal P}_{+} = \{\alpha \in \mathcal P; \alpha(v) > 0\}$. We can always find a subset ${\mathcal Q} \subset {\cal P}_{+}$ so that ${\mathcal Q}$ forms a basis for $\mathrm{t}^{*}$ and every $\alpha \in {\cal P}$ is an linear combination of elements of ${\mathcal Q}$ with integer coefficients, all of which are either nonnegative or nonpositive. The elements of ${\mathcal Q}$ are called {\it fundamental roots.}

It can be shown that every weight of $\pi$ is of the form $$\mu_{\pi} = \lambda_{\pi} - \sum_{\alpha \in \mathcal Q}n_{\alpha}\alpha$$
where each $n_{\alpha}$ is a non-negative integer and $\lambda_{\pi}$ is a weight of $\pi$ called the {\it highest weight}. Indeed if $\mu_{\pi}$ is any other weight of $\pi$ then $|\mu_{\pi}| \leq |\lambda_{\pi}|$. The highest weight of a representation is invariant under unitary conjugation of the latter and so there is a one-to-one correspondence between $\G$ and the space of highest weights $D$ of all irreducible representation of $G$. We can thus parameterise  $\G$ by $D$ and this a key step for Fourier analysis on nonabelian compact Lie groups. In fact $D$ can be given a nice geometrical description as the intersection of the weight lattice with the dominant Weyl chamber, but in order to save space we won't pursue that line of reasoning here. From now on we will use the notation $d_{\lambda}$ interchangeably with $d_{\pi}$ to denote the dimension of the space $V_{\pi}$ where $\pi \in \G$ has highest weight $\lambda$. For a more comprehensive discussion of roots and weights, see e.g. \cite{Fe} and \cite{Simon}.

\subsection{Sugiura Theory}

The main result of this subsection is Theorem \ref{Sug} which is proved in \cite{Sug}.

Let $M_{n}(\C)$ denote the space of all $n \times n$ matrices with complex entries and  ${\cal M}(G): = \bigcup_{\lambda \in D}M_{d(\lambda)}(\C)$. We define the {\it Sugiura space of rapid decrease} to be ${\cal S}(D): = \{F:D \rightarrow {\cal M}(G)\}$ such that
\begin{enumerate}
\item[(i)] $F(\lambda) \in M_{d(\lambda)}(\C)$ for all $\lambda \in D$, \item[(ii)]
$\lim_{|\lambda| \rightarrow \infty}|\lambda|^{k}|||F(\lambda)|||
= 0$ for all $k \in \mathbb{N}$. \end{enumerate}

${\cal S}(D)$ is a locally convex topological vector space with respect to the seminorms
$||F||_{s} = \sup_{\lambda \in D}|\lambda|^{s}|||F(\lambda)|||,$ where $s \geq 0$. We also note that $C^{\infty}(G)$ is a locally convex topological vector space with respect to the seminorms $||f|_{U} = \sup_{\sigma \in G}|Uf(\sigma)|$ where $U \in {\cal U}(\g)$, which is the universal embedding algebra of $g$ acting on $C^{\infty}(G)$ as polynomials in left-invariant vector fields on $G$, as described by the celebrated Poincar\'{e}-Birkhoff-Witt theorem.

\begin{theorem} \label{Sug} [Sugiura] There is a topological isomorphism between
$C^{\infty}(G)$ and ${\cal S}(D)$.
\end{theorem}

We list three useful facts that we will need in the next section. All can be found in \cite{Sug}.
\begin{itemize}
\item {\it Weyl's dimension formula} states that
$$ d_{\lambda} = \frac{\prod_{\alpha \in {\mathcal P}_{+}}(\lambda + \rho, \alpha)}{\prod_{\alpha \in {\mathcal P}_{+}}(\rho, \alpha)},$$
where $\rho:=\frac{1}{2}\sum_{\alpha \in {\mathcal P}_{+}}$ is the celebrated ``half-sum of positive roots''. From here we can deduce a highly useful inequality. Namely there exists $N > 0$ such that
\begin{equation} \label{WDineq}
             d_{\lambda} \leq N |\lambda|^{m}
\end{equation}
             where $m: = \#{\mathcal P}_{+} = \frac{1}{2}(d-r)$.

\item {\it Sugiura's zeta function} is defined by
$$ \zeta(s) = \sum_{\lambda \in D-\{0\}}\frac{1}{|\lambda|^{s}}$$
and it converges if $s > r$.

\item Let $(X_{1}, \ldots, X_{d})$ be a basis for $\G$ and let $\Delta \in {\cal U}(\g)$ be the usual Laplacian on $G$ so that
$$ \Delta = \sum_{i,j=1}^{d}g^{ij}X_{i}X_{j}$$
where $(g^{ij})$ is the inverse of the matrix whose $(i,j)$th component is $(X_{i}, X_{j})$. We may consider $\Delta$ as a linear operator on $L^{2}(G)$ with domain $C^{\infty}(G)$. It is essentially self-adjoint and $$ \Delta \pi_{ij} = - \kappa_{\pi} \pi_{ij}$$ for all $1 \leq i,j \leq d_{\pi}, \pi \in \G$, where $\pi \neq \delta \Rightarrow \kappa_{\pi} > 0$. The numbers $(\kappa_{\pi}, \pi \in \G\}$ are called the {\it Casimir spectrum} and if $\lambda_{\pi}$ is the highest weight corresponding to $\pi \in \G$ then
$$\kappa_{\pi} = ( \lambda_{\pi}, \lambda_{\pi} + 2\rho ).
$$From here we deduce that there exists $C > 0$ such that

\begin{equation} \label{Casin}
|\lambda_{\pi}|^{2} \leq \kappa_{\pi} \leq C(1 +
|\lambda_{\pi}|^{2}). \end{equation}

\end{itemize}

\subsection{Smoothness of Densities}

We can now establish our main theorem.

\begin{theorem} \label{smoden} $\mu \in {\cal P}(G)$ has a $C^{\infty}$ density if and
only if $\widehat{\mu} \in {\cal S}(D)$.
\end{theorem}

{\it Proof.} Necessity is obvious. For sufficiency its enough to show
$\mu$ has an $L^{2}$-density. Choose $s > r$ so that Suguira's zeta function converges. Then using Theorem \ref{sqrint} and (\ref{WDineq}) we have \bean
\sum_{\lambda \in
D-\{0\}}d_{\lambda}|||\widehat{\mu}_{\lambda}|||^{2} & \leq &
N\sum_{\lambda \in
D-\{0\}}|\lambda|^{m}|||\widehat{\mu}_{\lambda}|||^{2}\\
& \leq & N\sup_{\lambda \in
D-\{0\}}|\lambda|^{m+s}|||\widehat{\mu}_{\lambda}|||^{2}\sum_{\lambda
\in D-\{0\}}\frac{1}{|\lambda|^{s}}\\ & < & \infty.~~~~~~\Box \eean

We now investigate some classes of examples. We say that $\mu \in {\cal P}(G)$ is {\it central} if for all $\sigma \in G$, $$ \mu(\sigma A \sigma^{-1}) = \mu(A).$$ 
By Schur's lemma $\mu$ is central if and only if for each $\pi \in \G$ there exists $c_{\pi} \in \C$ such that
$$ \widehat{\mu}(\pi) = c_{\pi}I_{\pi}.$$

Clearly $m$ is a central measure. A standard Gaussian measure on $G$ is central where we say that a measure $\mu$ on $G$ is a {\it standard Gaussian} if it can be realised as $\mu_{1}^{(B)}$ in the convolution semigroup $(\mu_{t}^{(B)}, t \geq 0)$ corresponding to Brownian motion on $G$ (i.e. the associated Markov semigroup of operators is generated by $\frac{1}{2}\sigma^{2}\Delta$ where $\sigma > 0$.) For a more general notion of Gaussianity see e.g. \cite{He1}, section 6.2.  To verify centrality, take Fourier transforms of the heat equation to obtain $\widehat{\mu}(\pi) = e^{-\frac{1}{2}\sigma^{2}\kappa_{\pi}}I_{\pi}$ for each $\pi \in \G$.

Following \cite{App3} we introduce a class of central probability measures on $G$ which we call the {\it $CID_{\R}(G)$ class} as they are {\it central} and are induced by {\it infinitely divisible} measures on $\R$. Let $\rho$ be a symmetric infinitely divisible probability measure on $\R$ so we have the L\'{e}vy-Khintchine formula
$$ \int_{\R}e^{iux}\rho(dx) = e^{-\eta(u)}~\mbox{for all}~u \in \R$$
$$ \mbox{where}~\eta(u) = \frac{1}{2}\sigma^{2}u^{2} + \int_{\R-\{0\}}(1 -
\cos(u))\nu(du),$$  with $\sigma \geq 0$ and $\nu$ a L\'{e}vy measure, i.e. $ \int_{\R-\{0\}}( 1 \wedge
u^{2})\nu(du) < \infty$ (see e.g. \cite{Sa}.) We say $\mu \in CID_{\R}(G)$ if there exists
$\eta$ as above such that $$ \widehat{\mu}(\pi) =
e^{-\eta(\kappa_{\pi}^{\frac{1}{2}})}I_{\pi}~\mbox{for each}~\pi \in \G.$$ Examples of such measures are obtained by {\it subordination} \cite{Sa}. So let
$(\gamma_{t}^{f}, t \geq 0)$ be a subordinator with Bernstein
function $f$ so that for all $u \geq 0$
$$ \int_{0}^{\infty}e^{-us}\gamma_{t}^{f}(ds) = e^{-tf(u)}.$$
Let $(\mu_{t}^{(B)}, t \geq 0)$ be a Brownian convolution semigroup on $G$ (with $\sigma = \sqrt{2}$) so that for each $\pi \in \G$
$\widehat{\mu_{t}}(\pi) = e^{-t\kappa_{\pi}}I_{\pi}.$ then we obtain a
convolution semigroup of measures $(\mu_{t}^{f}, t \geq 0)$ in
$CIG_{\R}(G)$ by
$$ \mu_{t}^{f}(A) =
\int_{0}^{\infty}\mu_{s}^{(B)}(A)\gamma_{t}^{f}(ds)$$ for each $A \in {\cal B}(G)$ and we have
$$ \widehat{\mu_{t}^{f}}(\pi) = e^{-tf(\kappa(\pi))}I_{\pi}.$$
Examples (where we have taken $t=1$):
\begin{itemize}
\item Laplace Distribution $f(u) = \log(1 + \beta^{2}u)$,
$$ \widehat{\mu}(\pi) = (1 + \beta^{2}\kappa_{\pi})^{-1}I_{\pi}.$$

\item Stable-like distribution $f(u) =
b^{\alpha}u^{\frac{\alpha}{2}} (0 < \alpha < 2)$,
$$ \widehat{\mu}(\pi) =
e^{-b^{\alpha}\kappa_{\pi}^{\frac{\alpha}{2}}}I_{\pi}.$$

\end{itemize}

We now apply Theorem \ref{smoden} to present some examples of measures in the $CIG_{\R}$ class which have smooth densities (and one that doesn't).

\vspace{5pt}

Example 1. $\eta$ general with $\sigma \neq 0$ (i.e. non-vanishing Gaussian part)

\vspace{5pt}

Using (\ref{WDineq}) and (\ref{Casin}) we obtain \bean \lim_{|\lambda| \rightarrow \infty}|\lambda|^{k}|||\widehat{\mu}(\lambda)||| & = &
 \lim_{|\lambda| \rightarrow \infty}|\lambda|^{k}e^{-\eta{\kappa_{\pi}^{\frac{1}{2}}}}d_{\lambda}^{\frac{1}{2}}\\
 & \leq & \lim_{|\lambda| \rightarrow \infty}|\lambda|^{k}e^{-\frac{\sigma^{2}}{2}\kappa_{\lambda}}d_{\lambda}^{\frac{1}{2}}\\
 & \leq & N^{\frac{1}{2}}\lim_{|\lambda| \rightarrow \infty}|\lambda|^{k + \frac{m}{2}}e^{-\frac{\sigma^{2}}{2}|\lambda|^{2}} = 0.\eean

 \vspace{5pt}

 Example 2. Stable like laws are all $C^{\infty}$ by a similar argument.

 \vspace{5pt}

 Example 3. The Laplace distribution is not $C^{\infty}$. But it is $L^{2}$ if $r=1$ (e.g. $SO(3), SU(2), Sp(1)$.)

 \section{Deconvolution Density Estimation}

 We begin by reviewing the work of Kim and Richards in \cite{KR}. Let $X, Y$ and $\epsilon$ be $G$-valued random variables with $ Y = X\epsilon.$ Here we interpret $X$ as a signal, $Y$ as the observations and $\epsilon$ as the noise which is independent of $X$. If all three random variables have densities, then with an obvious notation we have $ f_{Y} = f_{X} * f_{\epsilon}.$ The statistical problem of interest is to estimate $f_{X}$ based on i.i.d. observations $Y_{1}, \ldots, Y_{n}$ of the random variable $Y$. We assume that the matrix $\widehat{f_{\epsilon}}(\pi)$ is invertible for all $\pi \in \G$. Our key tool is the {\it empirical characteristic function} $\widehat{f_{Y}}^{(n)}(\pi) :=\frac{1}{n}\sum_{i=1}^{n}\pi(Y_{i}^{-1})$. We then define the {\it non-parametric density estimator} (with smoothing parameters $T_{n} \rightarrow
\infty$ as $n \rightarrow \infty$) for $\sigma \in G, \nN$:
$$ f_{X}^{(n)}(\sigma): = \sum_{\pi \in \G:\kappa_{\pi} < T_{n}}d_{\pi}\tr(\pi(\sigma)\widehat{f_{Y}}^{(n)}(\pi)\widehat{f_{\epsilon}}(\pi)^{-1}).$$

The noise $\epsilon$ is said to be {\it super-smooth} of order $\beta > 0$ if there exists $\gamma > 0$ and $a_{1}, a_{2} \geq 0$ such that
$$||\widehat{f_{\epsilon}}(\pi))^{-1}||_{\infty} = O(\kappa_{\pi}^{-a_{1}}\exp(\gamma \kappa_{\pi}^{\beta}))~\mbox{and}~||\widehat{f_{\epsilon}}(\pi)||_{\infty} = O(\kappa_{\pi}^{a_{2}}\exp(-\gamma \kappa_{\pi}^{\beta})$$
as $\kappa_{\pi} \rightarrow \infty$. For example a standard Gaussian is super-smooth with $a_{i} = 0~(i = 1,2$). For $p > 0$, the Sobolev space
${\cal H}_{p}(G):=\{f \in L^{2}(G); ||f||_{p} < \infty\}$ where
$||f||_{p}^{2} = \sum_{\pi \in \widehat{G}}d_{\pi}(1 +
\kappa_{\pi})^{p}|||\widehat{f}(\pi)|||^{2}.$

 \begin{theorem}[Kim, Richards] If $f_{\epsilon}$ super-smooth of order $\beta$ and $||f_{X}||_{H_{s}(G)} \leq K$  for some $s > \frac{d}{2}$ where $K > 1$ then the optimal rate of convergence of $f_{X}^{(n)}$ to $f_{X}$ is $(\log(n))^{-\frac{s}{2\beta}}$.
 \end{theorem}

A natural question to ask is ``how smooth is super-smooth?" and we answer this as follows:

\begin{prop} If $f$ is super-smooth then it is smooth.
\end{prop}

{\it Proof.} For sufficiently large $\kappa_{\pi}$ and using (\ref{WDineq}) and (\ref{Casin}) we find that there exists $C > 0$ such that \bean |||\widehat{f}(\pi)||| & \leq & ||\widehat{f}(\pi)||_{\infty}|||I_{\pi}|||\\
& = & d_{\pi}^{\frac{1}{2}}||\widehat{f}(\pi)||_{\infty}\\
& \leq & N^{\frac{1}{2}}|\lambda_{\pi}|^{\frac{m}{2}}.C \kappa_{\pi}^{a_{2}}\exp(-\gamma \kappa_{\pi}^{\beta})\\
& \leq & K|\lambda_{\pi}|^{\frac{m}{2}}(1 + |\lambda_{\pi}|^{2})^{a_{2}}\exp(-\gamma |\lambda_{\pi}|^{2\beta}) \eean

from which it follows that $\widehat{f} \in {\cal S}(D)$ and the result follows by Theorem \ref{smoden}. $\hfill \Box$.

\end{document}